\documentclass{amsart}
\usepackage{hyperref}
\usepackage{txfonts}
\usepackage{amsfonts}

\newtheorem{theorem}{Theorem}[section]

\theoremstyle{definition}

\theoremstyle{remark}
\newtheorem{remark}[theorem]{Remark}

\numberwithin{equation}{section}

\begin{document}

\title{An extension of Bonnet-Myers theorem}

\author{}
\address{}
\curraddr{}
\email{}
\thanks{}

\author{Jianming Wan}
\address{School of Mathematics, Northwest University, Xi'an 710127, China}
\email{wanj\_m@aliyun.com}
\thanks{}


\date{}

\dedicatory{}


\begin{abstract}
We give a complementary generalization of the extensions of Bonnet-Myers theorem obtained by Calabi and also Cheeger-Gromov-Taylor.
\end{abstract}

\maketitle



The Bonnet-Myers states that a complete Riemannian manifold with  Ricci curvature $Ric_{M}\geq \delta>0$ is compact. In \cite{[C]} Calabi extended this by proving that, if for some point $p\in M$ every geodesic starting from $p$ has the property that $$\limsup_{a\rightarrow \infty}\int_{0}^{a} Ric_{M}(s)^{\frac{1}{2}}ds-\frac{1}{2}\ln(a)=\infty,$$ then $M$ is compact. In particular, this implies that (\cite{[WSY]} page 137) $M$ is compact provided  \begin{equation}\label{f0.1}
Ric_{M}(x)\geq\frac{1}{(4-\epsilon)(1+r)^{2}}
\end{equation}
 for $d(p,x)=r$ and all $r\geq 0$.
Cheeger, Gromov and Taylor (c.f. \cite{[CGT]} theorem 4.8) also proved a similar result in the same spirit. The basic idea of their proof is to study carefully the index form (or the second variation). The condition on Ricci curvature insures that the index form is negative in some direction. Then there exists conjugate points along any geodesic and $M$ has to be compact.

In this short note we give a complementary extension of Calabi and Cheeger-Gromov-Taylor's results. The main theorem is

\begin{theorem}
Let $M^{n}$ be a complete Riemannian manifold. If there exists $p\in M, k\geq2$ and $r_{0}>0$ such that
\begin{equation}\label{f0.2}
Ric_{M}(x)\geq\frac{C(n,k, r_{0})}{(r_{0}+r)^{k}}
\end{equation}
 for all $r\geq0$, where $d(p,x)=r$ and $C(n,k, r_{0})$ is a constant depending on $n,k, r_{0}$,
then $M^{n}$ is compact. In our situation, $C(n,k, r_{0})$ can be chosen to equal to $(n-1)\cdot\frac{(k-1)^{k}}{(k-2)^{k-2}}\cdot r_{0}^{k-2}$ for $k>2$ and
$(n-1)(1+\frac{r_{0}}{\epsilon}),\epsilon>0$ for $k=2$.
\end{theorem}

Since the case $ k<2$ is covered by \ref{f0.1}, we only consider $k\geq2$. It is easy to see that \ref{f0.1} or \ref{f0.2} covers the classical Bonnet-Myers theorem. If $Ric_{M}\geq \delta>0$, we can rescale the metric such that $\delta$ is bigger than the right hand of \ref{f0.1} or \ref{f0.2}.

From the viewpoint of index form or the second variation, the Ricci curvature decay rate 2 in \ref{f0.1} is the best possible. To guarantee that any geodesic from $p$ encounters conjugate points this is necessary. But to show that a complete Riemannian manifold is compact, ``meeting conjugate points'' is not need. Showing that the manifold has no ray is enough! This is our starting point. We would make use of \ref{f0.2} to show that $M^{n}$ contains no ray.

\section{A proof of the main theorem}
Assume that $M^{n}$ is noncompact. Then for any $p\in M$ there is a ray $\sigma(t)$ issuing from $p$.

Let $r(x)=d(p, x)$ be the distance function from $p$. We denote $A=Hess (r)$ outside the cut locus and  write $A(t)=A(\sigma(t))$.
The Riccati equation is given by
\begin{equation}\label{f1.1}
A^{'}+A^{2}+R=0.
\end{equation}
The A(t) is smooth except at $t=0$. Taking the trace we have
\begin{equation}\label{f1.2}
trA^{'}(s)+\|A(s)\|^{2}+Ric(T)=0,
\end{equation}
where $T=\sigma^{'}(s)$.

Integrate \ref{f1.2} over the interval $[\epsilon, t]$,

\begin{eqnarray*}
\int^{t}_{\epsilon}\|A(s)\|^{2}ds &=& trA(\epsilon)-trA(t)-\int^{t}_{\epsilon}Ric_{M}(T)ds\\
&< & \frac{n-1}{\epsilon}-\int^{t}_{\epsilon}\frac{C}{(r_{0}+s)^{k}}ds\\
&=& \frac{n-1}{\epsilon}-\frac{C}{k-1}[\frac{1}{(r_{0}+\epsilon)^{k-1}}-\frac{1}{(r_{0}+t)^{k-1}}].
\end{eqnarray*}
The ``$<$" holds from  $0\leq trA(t)< \frac{n-1}{t}$ and condition \ref{f0.2}. We claim that $0\leq trA(t)< \frac{n-1}{t}$.
Since $Ric_{M}>0 $, $trA(t)<\frac{n-1}{t}$ holds. To see $trA(t)\geq 0$, one can consider the excess function
$$e(x)=d(p,x)+d(\sigma(i), x)-i.$$
$e(x)\geq 0$ and $e(\sigma(t))\equiv 0$ for $0\leq t\leq i$. So $$\Delta e(\sigma(t))=\Delta (d(p,x)+d(\sigma(i), x))|_{\sigma(t)}\geq0.$$
We have $$trA(t)=\Delta d(\sigma(0),x)|_{\sigma(t)}\geq -\Delta d(\sigma(i), x))|_{\sigma(t)}\geq-\frac{n-1}{i-t}.$$
Let $i\rightarrow +\infty$. $trA(t)\geq 0$.

Let $t\rightarrow \infty$. The above integral inequality becomes
$$0\leq \int^{\infty}_{\epsilon}\|A(s)\|^{2}ds<\frac{(n-1)}{\epsilon}-\frac{C}{(k-1)(r_{0}+\epsilon)^{k-1}}.$$
We observe that when $C$ is very large, it is a contradiction. So $M$ must be compact. Now we work out the constant $C$ we need.
Solving $$\frac{(n-1)}{\epsilon}-\frac{C}{(k-1)(r_{0}+\epsilon)^{k-1}}\leq0,$$
we obtain
\begin{equation}\label{f1.3}
C\geq(n-1)(k-1)\frac{(r_{0}+\epsilon)^{k-1}}{\epsilon}.
\end{equation}
When $\epsilon=\frac{r_{0}}{k-2},k>2$, $\frac{(r_{0}+\epsilon)^{k-1}}{\epsilon}$ achieves minimal value. Substituting it into \ref{f1.3},  we have
$$C\geq(n-1)\cdot\frac{(k-1)^{k}}{(k-2)^{k-2}}\cdot r_{0}^{k-2}.$$
So we can choose $C(n,k,r_{0})=(n-1)\cdot\frac{(k-1)^{k}}{(k-2)^{k-2}}\cdot r_{0}^{k-2}$ for $k>2$. When $k=2$, $C(n,2,r_{0})=(n-1)(1+\frac{r_{0}}{\epsilon})$ for any $\epsilon>0$.

\begin{remark}
The estimate on the integral of $\|A(s)\|^{2}$ is inspired by Dai and Wei's paper \cite{[DW]}. In their work on Toponogov type comparison for Ricci curvature, the related estimate of Hessian plays an important role.
\end{remark}

\bibliographystyle{amsplain}

\begin{thebibliography}{10}
\bibitem {[C]} Calabi, E.,
\textit{On Ricci curvature and geodesics.} Duke Math. J. 34 1967 667-676.

\bibitem {[CGT]} Cheeger, J. ; Gromov, M. and Taylor, M., \textit{Finite propagation speed, kernel estimates for functions of the Laplace operator, and the geometry of complete Riemannian manifolds.} J. Differential Geom. 17 (1982), no. 1, 15-53.

\bibitem {[DW]} Dai, Xianzhe and Wei, Guofang \textit{A comparison-estimate of Toponogov type for Ricci curvature.} Math. Ann. 303 (1995), no. 2, 297-306.

\bibitem {[WSY]} Wu,H.; Shen, C. and Yu,Y., \textit{An introduction to Riemannian geometry (in Chinese).} Beijing University Press 1989.
\end{thebibliography}

\end{document}